\documentstyle[12pt]{article}
\def \C{{\Bbb C}}

\def \calB{{\cal B}}

\def \N{{\cal N}}

\def \O{{\Omega}}

\def \l{{\lambda}}
\def \Z{{\Bbb Z}}
\def \e{{\varepsilon}}
\def \tag{{\tilde{\frak{g}}}}
\def \ag{{\frak{g}}}

\def \o{{\omega}}
\def \K{{K\"{a}hler\ }}

\def \CC{{\cal C}}
\def \fd{{\bullet}}

\def \MN{{{\cal M}_{\cal N}}}
\def \s{{\sigma}}
\def \t{{\theta}}
\def \i{{\sqrt{-1}}}
\def \m{{\frak m}}
\def \d{{\partial}}

\def \tom{{\tilde \Omega}}
\def \tr{{\tilde{\cal R}}}
\def \tra{{\tilde{\cal R}_A}}

\def \tfa{{\tilde{\cal F}_A}}

\def \td{{\tilde D}}
\def \tc{{\tilde C}}
\def \V{{\bar V}}
\def \G{{\cal G}}
\def \g{{\gamma}}
\def \d{{\partial}}
\def \z{{\bar z}}
\def \a{{\alpha}}
\def \D{{\Delta}}
\def \b{{\delta}}

\def \proof{{\noindent{\it Proof.\ \ }}}
\newtheorem{Th}{THEOREM}[section]
\newtheorem{definition}[Th]{DEFINITION}
\newtheorem{prop}[Th]{PROPOSITION}

\newtheorem{lem}[Th]{LEMMA}

\title{Deformations of representations of fundamental groups of open K\"{a}hler
manifolds} \author{Philip A. Foth ${}^1$}
\begin{document}

\maketitle
\input amssym.def

\begin{abstract} Given a compact \K manifold, we consider the complement $U$ of
a divisor with normal crossings. We study the variety of unitary representations
of $\pi_1(U)$ with certain restrictions related to the divisor. We show that the
possible singularities of this variety as well as of the corresponding moduli
space of irreducible representations are quadratic. In the course of our proof
we exhibit a differential graded Lie algebra
which reflects our deformation problem. \end{abstract}

\footnotetext[1]{Sloan Doctoral Dissertation Fellow}

\tableofcontents

\section{Introduction}
\setcounter{equation}{0}

In the present paper we study the spaces of unitary representations of
fundamental groups of open \K manifolds. More precisely, let $X$ be a compact
\K manifold and let $U$ be the complement of a divisor with normal crossings
$D\subset X$, such that there exists a unique decomposition $D=\cup_{i=1}^r D_i$
into the union of $r$ smooth irreducible complex-analytic subvarieties of $X$.
In fact the Hironaka resolution of singularities theorem \cite{Hir}
states that every smooth quasi-projective variety is birational to a
manifold of this type.

Let us fix an $r$-tuple of conjugacy classes $\N=(\CC_1, ..., \CC_r)$ in
$U(N)$, a base point $b\in U$, simple loops $\g_i$ based at $b$ and encircling
$D_i$ for $1\le i\le r$ defining classes $[\g_i]\in\pi_1(U,b)$. We consider the
space $\tr(\pi_1(U,b), U(N))_\N$. of representations $\rho$ of $\pi_1(U, b)$
into $U(N)$ such that for $1\le i\le r$ the image $\rho([\g_i])$ belongs to the
conjugacy class $\CC_i$. This space has the natural structure of a real
algebraic variety.

For an algebraic variety $Y$ and a point $y\in Y$ wee say that $Y$ is quadratic
at $y$ if there exists an analytic embedding $\phi$ of a neighborhood $U$ of $y$
into a vector space such that $\phi(y)=0$ and $\phi(U)$ is given as the zero
locus of a finite number of homogeneous quadratic equations. We prove

\

{\bf THEOREM \ref{Th:mt} } {\it Let $\rho: \pi_1(U,b)\to U(N)$ be a
representation such that $\rho([\g_i])\in\CC_i$ for $1\le i\le r$, then the
space $\tr(\pi_1(U,b), U(N))_\N$ is quadratic at $\rho$.}

\

This theorem was proved by Goldman and Millson \cite{GM2} in the compact case,
i.e.  when $D$ is empty. Working in the context of integrable connections,
Biquard \cite{Biq} showed that when $r=1$ (i.e. smooth divisor case) the
singularities of the space of isomonodromic deformations of an integrable
logarithmic connection are quadratic.

If we restrict our attention to irreducible representations, then we can
consider the moduli space $$\MN:= \tr^{irr}(\pi_1(U,b), U(N))_\N/U(N),$$ where
the group $U(N)$ acts on the space of representations by conjugation. This space
was the subject of a work of Jean-Luc Brylinski and the author \cite{BF}. In
fact, it was shown there that when a certain obstruction vanishes, the tangent
space to $\MN$ at a point $[\rho]\in\MN$ is identified with the first
intersection cohomology group $IH^1(X, \tag)$, where $\tag$ is the local system
on $U$ corresponding to the representation $Ad\circ\rho : \pi_1(U,b)\to {\frak
u}(N)$. This was used to exhibit a symplectic form on the smooth locus of
$\MN$, which turns out to be a \K form when $X$ is projective. Another possible
description of the tangent space $T_{[\rho]}\MN$ is the first relative group
cohomology group $H^1(\pi_1(U), (\Gamma_i), {\frak u}(N))$, where $\Gamma_i$ be
a cyclic subgroup of $\pi_1(U)$ generated by $\g_i$ and $Ad\circ\rho$ give
${\frak u}(N)$ the structure of a $\pi_1(U,b)$-module.

A straightforward consequence of the above Theorem is that if $\rho$ is an
irreducible representation satisfying the above hypotheses, then $\MN$ is
quadratic at $[\rho]$.

Let us take a representation $\rho$ defining a point $[\rho]\in\MN$. Consider
the pairing $$ H^1(\pi_1(U), (\Gamma_i), {\frak u}(N))\times H^1(\pi_1(U),
(\Gamma_i), {\frak u}(N)) \to H^2(\pi_1(U), (\Gamma_i), {\frak u}(N))$$ given
by the cup product in relative group cohomology together with the Lie bracket as
a coefficient pairing.

\

{\bf THEOREM \ref{Th:pair} } (cf. \cite{GM2}) {\it A point $[\rho]\in\MN$ is a
smooth point of $\MN$ if the above pairing is identically zero.}

\

Following the ideology developed by Schlessinger-Stasheff \cite{SS},
Deligne, Goldman-Millson \cite{GM2} we construct a differential graded Lie
algebra (DGLA) $B^{\fd}(\tag)$ which is the controlling DGLA for our
deformation problem. This algebra consists of smooth
differential forms on $U$ with coefficients in $\tag$ with certain conditions
on asymptotics of the forms along $D$.
We say that a DGLA is formal if it is quasi-isomorphic to its
cohomology algebra.

\

{\bf THEOREM \ref{Th:form} } {\it The DGLA $B^{\fd}(\tag)$ is formal. }

\

This is similar to the fundamental result of
Deligne-Griffiths-Morgan-Sullivan \cite{DGMS} about the formality of the de Rham
algebra of a compact \K manifold. One of the important ingredients in the
course of our proof is an analogue of the $\bar\d$-Poincar\'e lemma. The DGLA
$B^{\fd}(\tag)$ is naturally bigraded by holomorphic and anti-holomorphic
degrees. Our Lemma \ref{lem:d'l} tells us that locally the complex $$\cdots
\stackrel{\nabla''}{\to}B^{p,q-1}(\tag)\stackrel{\nabla''}{\to}
B^{p,q}(\tag)\stackrel{\nabla''}{\to}B^{p,q+1}(\tag)\stackrel{\nabla''}{\to}
\cdots$$ is exact.

The section of our paper which relates this DGLA and our deformation problem
heavily relies on the theory which Goldman and Millson developed in the compact
case. The deformation theory that we study is equivalent to the deformation
theory of flat unitary connections such that the monodromy around $i$-th
irreducible component of the divisor lies in a specified conjugacy class
$\CC_i$. As in \cite{GM2} we formulate the equivalence theorem in terms of
functors from the category of Artin local algebras.

We must mention that Kapovich and Millson in \cite{KM} studied the relative
deformation theory of representation with relation to the local deformations of
linkages.

\

{\bf Acknowledgments.} I am very grateful to Jean-Luc Brylinski for the
excellent guidance which I and this paper indubitably benefited from.
I would like to thank John Millson for important conversations, James Stasheff
for useful comments, and Indranil Biswas for drawing my attention to the results
of Timmerscheidt.

\section{Forms with logarithmic poles}
\setcounter{equation}{0}

Let $X$ be a compact \K manifold of complex dimension $d$ endowed with a \K
form $\l$ and let $D$ be a normal crossing divisor on $X$ which can be
represented as a union of $r$ smooth irreducible complex-analytic subvarieties
of codimension $1$:  $$D=\bigcup_{i=1}^rD_i.$$ We let $G=U(N)$ and $\ag={\frak
u}(N)$ and we pick a non-degenerate symmetric invariant bilinear form $B\langle,
\rangle$ on $\ag$.  We also fix a set of $r$ conjugacy classes in $G$: $${\cal
N}=(\CC_1, ..., \CC_r),$$ one for each component of the divisor. We let
$U:=X\setminus D$, we let the map $j: U\hookrightarrow X$ be the inclusion, and
we pick a base point $b\in U$ and we denote $\pi_1(U):=\pi_1(U, b)$. Let $\g_i$
be a simple loop encircling $D_i$ for $1\le i\le r$ and let $\Gamma_i\simeq\Z$
be a subgroup of $\pi_1(U)$ generated by $[\g_i]$ - the class of $\g_i$ in
$\pi_1(U)$.

We consider the space $\tr(\pi_1(U), G)_\N$  of representations $$\rho:
\pi_1(U)\to G$$ such that $\rho([\g_i])\in\CC_i$. This space has the structure
of a real algebraic variety. We also consider the moduli space $\MN$ of
irreducible representations of this type, i.e. $$\MN=\tr^{irr}(\pi_1(U),
G)_\N/G,$$ where $G$ acts by conjugation. In general these spaces are not
smooth, and we shall show that in fact all possible singularities are quadratic:
i.e.  a neighbourhood of each point $[\rho]\in\MN$ is analytically isomorphic to
a neighbourhood of the origin of ${\cal L}$ in $E$. Here $E$ is a vector space
and ${\cal L}$ is the zero locus of a quadratic form $Q$ in $E$ which in general
is given by several quadratic equations. A point is smooth then if $Q=0$
identically on $E$.

Let $D^{(m)}$ be the subvariety consisting of points which belong to at least
$m$ irreducible components of the divisor $D$, e. g. $D^{(0)}=X$ and
$D^{(1)}=D$.  Let $\td^{(m)}\to D^{(m)}$ be normalization maps and let
$$v^m:\td^{(m)}\to D^{(m)}\hookrightarrow X$$ be th composition of
normalizations with inclusions.  We also let $\tc^{(m)}:=(v^m)^{-1}(D^{(m+1)})$,
and it is a divisor with normal crossings on $\td^{(m)}$ (possibly empty).

Given a representation $\rho$ such that $[\rho]\in\MN$ we can consider $\ag$ as
a $\pi_1(U)$-module via the adjoint representation followed by $\rho$. We also
let $\tag$ stand for the local system corresponding to $\ag$, let $V$ be
the flat bundle over $U$ corresponding to the local system $\tag$, and finally
let $(\V , \nabla)$ be the Deligne extension of $V$. Here $\bar V$ is a
holomorphic vector bundle over $X$ and $\nabla$ is a holomorphic connection with
logarithmic singularities along $D$. We refer the reader to \cite{De1} for the
details. The connections are extended as usual to holomorphic forms with
logarithmic poles.

There exist \cite{DelHod} higher residues $Res_m(\V )$ which give
morphisms of complexes: $$Res_m(\V ): \ \ \O^{\fd}_X\langle D\rangle
\otimes\V \to v^m_*(\O^{\fd}_{\td^{(m)}}\langle \tc^{(m)}\rangle\otimes {\bar
V}_m)[-m].$$ Here $\V _m$ is the unique vector subbundle of $(v^m)^*\V $
equipped with a unique holomorphic integrable connection $\nabla_m$ with
logarithmic poles along $\tc^{(m)}$ such that
$$Ker{\nabla_m}_{|\td^{(m)}\setminus
\tc^{(m)}}=(v_m)^{-1}((j_*\tag)_{|\td^{(m)}\setminus\td^{(m+1)}}).$$ This means
that $$(\V _m, \nabla_m)$$ is the canonical extension of
$$(v_m)^{-1}((j_*\tag)_{|\td^{(m)}\setminus\td^{(m+1)}}).$$ As usual, we denote
by $\O^{\fd}_X\langle D\rangle$ the complex of sheaves of holomorphic forms on
$U$ with logarithmic poles along $D$. The fact that $Res_m(\V )$ is a
morphism of complexes means that $$Res_m(\V )\circ\nabla=\nabla_m\circ Res_m(\V
).$$ Let us define $$\tom^{\fd}_X(\V ):=Ker Res_1({\bar V});$$ it is a complex
of sheaves of holomorphic forms on $X$ with coefficients in $\V $ which have
logarithmic poles along $D$ and have no residues. Residue maps take values in
the monodromy invariant part of the local system, so when all the eigenvalues of
monodromy transformation are different from $1$, the complexes
$\tom^{\fd}_X(\V)$ and $\O^{\fd}_X\langle D\rangle\otimes \V$ coincide.

The complex $\tom^{\fd}_X(\V )$ admits a different description \cite{T'}.
First, on $U$ we let $\tom^p_X(\V)_{|U}:=\O^p_U\otimes\V_{|U}$. For $x\in D$ we
pick a small polycylinder $\D$ containing $x$ such that $D\cap\D$ is given by
$z_1\cdots z_l=0$ in local holomorphic coordinates $z_1, ..., z_d$. Let $T_i$,
$1\le i\le l$, be the monodromy transformation of $\tag_{|\D\setminus D}$
around $D_i$. Let $\tom^1_X(\V)_{|\D}\subset\O^1_X\langle
D\rangle\otimes\V_{|\D}$ be generated over ${\cal O}_{\D}$ by
$\O^1_{\D}\otimes_{{\cal O}_{\D}}\V_{|\D}$ and $$\O^1_{\D}\langle D_i\rangle
\otimes_{\C}(Ker(T_i-Id))^{\perp},$$ where $(Ker(T_i-Id))^{\perp}$ is the
orthogonal complement of the local subsystem $Ker (T_i-Id)$ of
${\tag}_{\D\setminus D}$, and $1\le i\le l$. Similarly we define the groups
$\tom^p_X(\V)_{|\D}$. Those sheaves glue together
nicely to a subcomplex $(\tom^{\fd}_X(\V), \nabla)$ of the complex of sheaves of
holomorphic differential forms on $X$ with coefficients in $\V$ with logarithmic
poles along $D$.

We will use the notation $d'$ and $d''$ respectively for the holomorphic and
anti-holomorphic covariant derivatives corresponding to $\nabla$.

We define ${\cal L}^{p,q}(\tag)_{(2)}$ to be
the sheaf of measurable $(p,q)$-forms $\a$ on $U$ with values in $\tag$ such
that both $\a$ and $(d'+d'')\a$ are locally square-integrable. Here as usual we
have the Poincar\'e metric near $D$ and the \K metric on $U$ which are
compatible. It means that the \K metric near $x\in D$ has the same asymptotic
form as the Poincar\'e metric on $\D\cap D$, where $\D$ is a small polycylinder
centered at $x$ (cf. \cite{Zu}). Set $$L^{p,q}(U, \tag)_{(2)}:=\Gamma(U, {\cal
L}^{p,q}(\tag)_{(2)}).$$ We let $$H^{p,q}(U,
\tag)_{(2)}:=H^q(L^{p,\fd}(U, \tag)_{(2)}, d'').$$ Analogously one defines the
space $H^k(U, \tag)_{(2)}$. Let us recall the following result due to
Timmerscheidt \cite{T'}:

\begin{prop} \noindent (a) $$H^{p,q}(U,\tag)_{(2)}$$ is conjugate-isomorphic to
$$H^{q,p}(U, \tag')_{(2)},$$ where $\tag'$ is the local system dual to $\tag$
and $$H^k(U,\tag)_{(2)}\simeq\bigoplus_{p+q=k}H^{p,q}(U, \tag)_{(2)}.$$

(b) The complexes of sheaves $\tom^p_X(\V)$ and ${\cal L}^{p,
\fd}(\tag)_{(2)}$ are quasi-isomorphic; thus $$H^q(X, \tom^p_X(\V))\simeq
H^{p,q}(U, \tag)_{(2)}.$$

(c) The spectral sequence $$E^{p,q}_1=H^q(X, \tom^p_X(\V))$$ degenerates at
$E_1$, abuts to $H^{p+q}(X, j_*\tag)$ with $$H^k(X, j_*\tag)\simeq\bigoplus
H^q(X, \tom^p_X(\V)).$$

(d) We have a conjugate linear isomorphism $$H^q(X, \tom^p_X(\V))\simeq H^p(X,
\tom^q_X(\bar{V'})), $$ where $\bar{V'}$ is the canonical extension of $V'$ -
the flat bundle corresponding to $\tag'$. \label{prop:Tim} \end{prop}

\section{Construction of the DGLA and formality}
\setcounter{equation}{0}

In this section everything is defined over the field of complex numbers.
A {\it differential graded Lie algebra} (DGLA) $(L, d)$ consists of a graded Lie
algebra $$L=\oplus_{i\ge 0}L^i, \ \ [ , ]: L^i\times L^j\to L^{i+j}, $$
satisfying for $\a\in L^i, \beta\in L^j$, and  $\g\in L^k$:
$$[\a, \beta]+(-1)^{ij}[\beta, \a]=0,$$
$$(-1)^{ki}[\a,[\beta,\g]]+(-1)^{ij}[\beta,[\g,\a]]+(-1)^{jk}[\g,
[\a,\beta]]=0$$ together with a derivation $d$ of degree $1$ (also called
a differential):  $$d: L^i\to L^{i+1}, \ d\cdot d=0, \ d[\a, \beta]=[d\a,
\beta]+(-1)^i[\a, d\beta].$$

The cohomology of any
DGLA is a DGLA too, considered with zero differential. We say
that two DGLAs $(L_1, d_1)$ and $(L_2, d_2)$ are {\it quasi-isomorphic} if there
exists a third DGLA $(L_3, d_3)$ and DGLA homomorphisms $i$ and $p$: $$(L_1,
d_1)\stackrel{i}{\leftarrow} (L_3, d_3) \stackrel{p}{\to} (L_2, d_2)$$ such that
both $i$ and $p$ induce isomorphisms in cohomology. We also say that a DGLA is
{\it formal} if it is quasi-isomorphic to its cohomology.

Next we construct a double complex version of a DGLA
$(B^{\fd,\fd}(\tag), d', d'')$ such that the corresponding ordinary DGLA
$(B^{\fd}(\tag), d'+d'')$ controls deformations of $[\rho_0]\in\MN$. The
space $B^{p,q}(\tag)$ is defined as the space of global sections of a sheaf
$\calB^{p,q}(\tag)$. If $x\in U$, and $\D$ is a small polydisc around $x$, then
the sections of $\calB^{p,q}(\tag)(\D)$ are the smooth differential forms of
type $(p,q)$ on $\D$ with coefficients in $\tag$. Let us now consider a point
$x\in D$ and a polydisc $\D$ containing this point such that $\D\cap D$ is given
in local coordinates $z_1, ..., z_d$ by $z_1\cdots z_k=0$. Locally, since we
have a semi-simple local system and a normal crossing divisor, the local
fundamental group is commutative (and isomorphic to $\Z^k$) and the local system
splits into the direct sum of local systems of rank one. Therefore, we may
assume that $\tag$ has rank $1$. We also may assume that the monodromy around
$D_i:=\{ z_i=0\}$ is non-trivial with the eigenvalue $e^{-\i\t_i}$ with $0< \t_i
<2\pi$.  Otherwise, we can ignore the components of the divisor with
trivial monodromy. Any ordinary homogeneous differential form $\o$ of type
$(p,q)$ we represent as $$\o=f(z_1, ..., z_d)dz_1^{\b_1}\wedge\cdots\wedge
dz_d^{\b_d}\wedge d\z_1^{\s_1}\wedge\cdots\wedge d\z_d^{\s_d},$$
$$\b_i, \s_j \in\{ 0,1\}, \ \ \ \sum \b_i=p, \ \ \sum \s_i=q.$$

\begin{definition} We say that a differential form $\o$ of type 
$(p,q)$ given by $$\o=f(z_1, ..., z_d)dz_1^{\b_1}\wedge\cdots\wedge 
dz_d^{\b_d}\wedge d\z_1^{\s_1}\wedge\cdots\wedge d\z_d^{\s_d}$$ satisfies {\rm 
condition C} if $f(z_1, ..., z_d)$ is smooth over $\D\setminus D$, is bounded 
for $p=q=0$, and for $p+q >0$ has the following asymptotics: $$f(z_1, ...,
z_n)=O(r_1^{\e-\b_1-\s_1}\cdots r_k^{\e-\b_k-\s_k})$$ for some $\e >0$.
\end{definition}

If we have locally on $\D\setminus D$ a multi-valued horizontal section $s$:
$\nabla s=0$ then the monodromy operator $T_i$ corresponding to $D_i$ acts on it
as $T_is=e^{-\i\t_i}s$. Now we will say that  $\a\otimes s\in
\calB^{p,q}(\tag)(\D)$ for a homogeneous differential form $\a$ of type $(p,q)$
if $\a$, $d'\a$, $d''\a$, and $d'd''\a$ all satisfy condition C.

It is clear that defined in such a way spaces of sections glue nicely together
to give a double complex of sheaves $\calB^{\fd, \fd}(\tag)$ and we get a double
complex $(B^{\fd, \fd}(\tag), d', d'')$ by taking the spaces of global sections 
of the above double complex of sheaves. One easily checks that the cup-product 
of differential forms combined with Lie bracket as the coefficient pairing
endowes $B^{\fd}(\tag)$ with the structure of DGLA.

Our next goal is to prove the following
\begin{prop} The spectral sequence defined by the filtration associated to
either degree of the double complex $(B^{p,q}(\tag), d', d'')$ degenerates at
$E_1$. The two induced filtrations on $H^k$ are $k$-opposite.
\label{prop:pss} \end{prop}

We need to establish certain results analogous to
the classical $d''$-Poincar\'e Lemma. First, we prove the result in
complex dimension one. Let $\D\subset\C$ be a a disc centered at zero,  let
${\bar\D}$ be its closure, and let $\D^*$ be the disc $\D$ punctured at zero.
The number $\e$ is always assumed to be a small positive real number.

\begin{lem} Given a complex $C^{\infty}$ function $f$ in an open punctured
neighbourhood $W_1^*$ of ${\bar\D}$ such that the asymptotics of $f$ at zero is
as $f=O(r^{\e-k})$, there exists a $C^{\infty}$ function $g=O(r^{\e-k+1})$ in an
open punctured neighbourhood $W_2^*\subset W_1^*$ of $\bar\D$ such that $${\d
g\over\d\z}=f$$ in $W_2^*$. Moreover, if $f$ is $C^{\infty}$ or has asymptotics
like above in some additional parameters, then $g$ can be chosen to have the
same properties.  \label{lem:l00} \end{lem}

\proof The main idea of the proof is that the controlled asymptotics at the
origin allows us to use Fourier analysis, which turns out to be an appropriate
tool in this context. We notice that the function $h(z)=r^{4-\e}f$ is of class
$C^2$, and $h$, ${\d h\over\d r}$, and ${\d^2 h\over \d r^2}$ vanish at the
origin.

We decompose $$h(z)=\sum_{n\in\Z} e^{\i n\phi}h_n(r),$$ where
$$h_n(z)={e^{-\i n\phi}\over2\pi}\int_0^{2\pi}h(e^{\i\phi}z)e^{-\i n\phi}
d\phi$$ (clearly, it is a function of $r$). Each
function $h_n(r)$ is also of class $C^2$ with vanishing $h_n$, ${dh_n\over d
r}$, and ${d^2 h_n\over d r^2}$ at the origin.  Moreover, integration by parts
shows that $$\sup ||h_n(r)||\le {C\over n^2}\sup ||{\d^2 h\over \d r^2}||$$
for a universal constant $C$, which
shows that the above series is absolutely convergent.

Now we let $$g_n(r)=-r^{4-\e-n}\int_r^1 2 \rho^{n-4+\e}h_n(\rho)d\rho$$ and
we immediately see that $g_n(r)$ is of class $C^3$ and vanishes at zero
together with its first, second, and third derivatives. Furthemore, one sees
easily that the series $$j(z)=\sum_{n\in\Z} e^{\i(n-1)\phi} g_n(r)$$ is
absolutely convergent together with its first derivative on compact subsets of
$W_1$. Finally, we let $g(z)=j(z)r^{\e-4}$ and using the identity
$${\d\over\d\z}= {e^{\i\phi}\over 2}{\d\over\d r}+{\i e^{\i\phi}\over
2r}{\d\over \d\phi}$$ we conclude that the function $g(z)$ is a well-defined
smooth $C^{\infty}$ function in $W_2^*$ and satisfies ${\d g(z)\over
\d\z}=f(z)$. The fact that the function $g(z)$ has the desired asymptotics is a
consequence of our explicit formulae. $\bigcirc$

\

The following statement is local in nature and is a direct consequence of the
above Lemma \ref{lem:l00} and the standard separation of variables arguments
which can be found e.g. in \cite{Cun} Vol.I E, Theorem 3.

\begin{lem} For $q>0$, $x\in X$ let $\o$ be an element of
$B^{p,q}(\tag)$ defined in an open neighbourhood of a polydisc $\bar \D$
containing $x$ such that $d''\o=0$.  Then there exists
$\a\in B^{p,q-1}(\tag)$ defined in an open neighbourhood of $\bar \D$ such
that $\o=d''\a$. \label{lem:d'l} \end{lem}

This Lemma allowes us to identify $E_1$ of the spectral sequence
of the double complex $B^{p,q}(\tag)$ corresponding to the holomorphic
filtration.  The local exactness of the columns of the
$E_0^{p,q}:=B^{p,q}(\tag)$ is equivalent to Lemma \ref{lem:d'l} and the
statement that the kernel
of the sheaf homomorphism $$d'': \calB^{p,0}(\tag)\to \calB^{p,1}(\tag)$$ is
the sheaf $\tom^p_X(\V)$. For simplicity, we will discuss the case $d=1$ in a
polydisc $\D\ni x\in D$ when the monodromy around $z:=z_1=0$ is equal to
$e^{-\i\t}$, where $0\le\t <2\pi$. Thus if we have a multi-valued horizontal
section $s$: $\nabla s=0$ then the monodromy operator $T$ acts on it as
$Ts=e^{-\i\t}s$. Let us consider the uni-valued section $\mu=z^{\t/2\pi}s$ which
spans the canonical (Deligne \cite{De1}) holomorphic extension $\V$ of $\tag$.
We have: $$\nabla\mu={\theta\over 2\pi}{dz\over z}\otimes\mu$$ and the Leibnitz
rule also imply that for any local section $f(z)\in\calB^{1,0}(\tag)(\D)$, we
have $$d''(f\otimes \mu)=d''f\otimes \mu$$ and the asymptotics of $f\otimes \mu$
is clearly $O(r^{\e-1-{\t\over 2\pi}})$. Thus $d''f=0$ implies that $f$ is
necessarily a holomorphic $1$-form with at most logarithmic singularities at
$z=0$ when $\t>0$, and is a holomorphic one-form in $\D$ when $\t=0$.

Similar arguments apply straightforwardly in higher dimension situations to
prove $$E_1^{p,q}=H^q(X, \tom^p_X(\V)).$$

The construction of our bigraded DGLA $B^{\fd, \fd}(\tag)$ was performed in
order to ensure a certain symmetry with respect to complex conjugation.  In
particular, if we consider our double complex with the anti-holomorphic
filtration, then the spectral sequence $\bar E$ associated to it will have
${\bar E}^{p,q}_1= H^p({\bar X}, {\tilde{\Omega}}^q_{\bar X}({\bar \V}))$.
Here ${\bar\V}$ is the anti-holomorphic extension of $\tag$ and
$\bar X$ is the manifold $X$ considered with the complex conjugate complex
structure. For example, when $d=1$ and we have rank one local system
$\tag$ with horizontal multi-valued section $s$ as above then $\bar\V$
is spanned by $\z^{1-{\t\over 2\pi}}s$.

It follows that each of our spectral sequences identifies with the spectral
sequence for the $L_2$ double complex, considered in Proposition
\ref{prop:Tim}. More precisely, $E_1^{p,q}$ identifies with $H^{p,q}(U,
\tag)_{(2)}$, and ${\bar E}^{p,q}_1$ with $H^{p,q}(U, \tag')_{(2)}$.
Hence, by Proposition \ref{prop:Tim} the spectral sequences degenerate at
$E_1$ and the induced filtrations on $E_1={\bar E}_1=E_{\infty}$ produce
$k$-opposite filtrations on $H^k(U, \tag)_{(2)}$. This completes the proof
of Proposition \ref{prop:pss} $\bigcirc$

\

Next, we apply Proposition 5.17 from \cite{DGMS} which tells us that
our Proposition \ref{prop:pss} is in fact equivalent to the $d'd''$-Lemma for
the DGLA $B^{p,q}(\tag)$. Repeating verbatim the proof of The Main Theorem
from the same source \cite{DGMS} then assures us that for the diagram
$$(B^{\fd}(\tag), d'+d'')\stackrel{i}{\leftarrow}(Ker
[d':B^{\fd}(\tag) \to B^{\fd}(\tag)], d'')\stackrel{p}{\to}(H^{\fd}(X,
j_*\tag), 0)$$ the maps $i^*$ and $p^*$ induce isomorphisms on cohomology. We
conclude with \begin{Th} The DGLA $(B^{\fd}(\tag), d'+d'')$ is formal.
\label{Th:form} \end{Th}

In fact, the complex of $L_2$ differential forms
$L^{\fd}(U, \tag)_{(2)}$ is not necessarily itself a DGLA, and this is the main
reason that we had to construct the DGLA $B^{\fd}(\tag)$ which has the right
cohomology and therefore adequately reflects our deformation problem.
We further notice that we have a natural inclusion of double complexes:
$B^{p,q}(\tag)\hookrightarrow L^{p,q}(U, \tag)_{(2)}$. We immediately have

\begin{prop} This inclusion induces an isomorphism
on cohomology. \end{prop}

\section{The deformation theory of flat bundles}
\setcounter{equation}{0}

In this section we follow ideas in the paper of Goldman and Millson
\cite{GM2}. Basically we shall adapt key results of their works to our
situation, which is different, because our \K manifold is not compact and
we have restrictions on monodromy transformations around $D_i$'s. We shall adopt
categorical language and our major claims will be stated as equivalences
of certain groupoids. For complete and exhaustive treatment applied in many
situations we refer the reader to \cite{GM2}.

First, we introduce an Artin local algebra $A$ over $\C$ (which has residue
field $\C$) with maximal ideal $\m$. We assume $A$ to be unital so there are
well-defined maps $$\C\stackrel{i}{\to}A\stackrel{p}{\to}{\C}$$ of inclusion and
projection onto the residue field.

Let us have a unitary irreducible representation $\rho_0: \pi_1(U)\to G$
defining a point $[\rho]\in\MN$, i.e. $\rho_0$ sends the class $[\g_0]$ of a
designated simple loop $\g_i$ around $D_i$ to the fixed conjugacy class
$\CC_i\subset G$ for $1\le i\le r$. This representation defines a flat unitary
bundle $(W, \nabla_0)$ on $U$ with the property that the monodromy
transformation around $D_i$ lies in $\CC_i$, $1\le i\le r$. In this situation
the bundle $V$ corresponds to the bundle of skew-adjoint endomorphisms of $W$.

Let us consider the groupoid $\tfa(\nabla_0)$ such that

\noindent - its objects are $A$-linear integrable smooth connections $\nabla$
over $U$ with values in $V\otimes A$

\noindent - these connections deform $\nabla_0$, i.e. $p(\nabla)=\nabla_0$

\noindent - any element $\nabla\in\tfa(\nabla_0)$ has the property that
$\a:=i(\nabla_0)-\nabla$, $d'\a$, $d''\a$, and $d'd''\a$ satisfy condition C 
from Section 3, where the map $i$ is induced by the inclusion $i: \C\to A$.

The morphisms of this groupoid are given by the group of gauge equivalences
$\G_A$ with values in $A$, i.e. smooth bundle automorphisms of $V\otimes A$ over
$X$ which preserve $\nabla$.

Next we recall the DGLA $(B^{\fd}(\tag), d'+d'')$ and the quadratic map
$$Q_A: B^{1}(\tag)\otimes\m\to B^{2}(\tag)\otimes\m, \ \ Q_A(\a)=(d'+d'')
\a+{1\over 2}[\a, \a].$$ The subspace $Q^{-1}_A(0)$ is preserved by the natural
action of $\exp(B^0(\tag)\otimes\m)$. This action makes the set of objects
$Q^{-1}_A(0)$ into a groupoid ${\tilde{\cal C}_A}$. We also notice that
$d'+d''$ is the decomposition of the covariant derivative corresponding to
$\nabla_0$.

Given any $\nabla_0$, Goldman and Millson in \cite{GM2} define uniquely
an element $\tilde{\nabla}_0\in\tfa(\nabla_0)$ called the {\it extension} of
$\nabla_0$. Take $\nabla\in\tfa(\nabla_0)$; the difference
$\nabla-\tilde{\nabla}_0$ is an element of $B^{1}(\tag)\otimes\m$, because it
clearly has the right asymptotics modulo $\m$ by definition of $\tfa(\nabla_0)$.
Proposition 6.6 from \cite{GM2} now translates as

\begin{prop} The above correspondence defines an isomorphism of groupoids
$$\tfa(\nabla_0)\to{\tilde{\cal C}_A}$$ depending naturally upon $A$. \end{prop}

To throw a bridge between deformations of a flat irreducible connection
$\nabla_0$ and the corresponding representation $\rho_0$ we shall need another
groupoid $\tra(\rho_0)$ whose objects are representations of $\pi_1(U,b)$ into
$G_A$ - the group of $A$-points of $G$ which deform $\rho_0$ (i.e.
$p(\rho)=\rho_0$) subject to the condition that for $\rho\in\tra(\rho_0)$ the
element $\rho(\g_i)$ is conjugate to $\rho_0(\g_i)$ for $1\le i \le r$.
We also require that if the intersection of several components $D_{i_1}, ...,
D_{i_s}$ of the divisor is non-empty, then $\rho(\g_{i_1}), ...,
\rho(\g_{i_s})$ are simultaneously conjugate in the group $G_A$ to
$\rho_0(\g_{i_1}), ..., \rho_0(\g_{i_s})$ respectively. The set of morphisms in
this groupoid is provided by the action of the group $\exp(\ag\otimes\m)$.

We shall need the following well-known direct consequence of Nakayama's lemma:

\begin{lem} Let $A$ be a noetherian local algebra over $\C$ and
let $f_A: M_A\to N_A$ be an $A$-linear map between two free finitely generated
$A$-modules. Then if $$f_A\otimes Id: M_A\otimes_A \C \to N_A \otimes_A\C $$ is
an isomorphism, then $f_A$ is an isomorphism as well. \label{lem:lNak} \end{lem}

\

Let us define the map $\e_b: \G_A \to G_A$ given as the evaluation at the base
point $b$. We recall the monodromy functor \cite{GM2} (where it is called
holonomy, though) $(mon_b, \e_b)$ between $\tfa(\nabla_0)$ and $\tra(\rho_0)$
defined naturally. Now we will adapt Proposition 6.3 from \cite{GM2} to our
purposes:

\begin{prop} The functor $$(mon_b, \e_b): \tfa(\nabla_0) \to \tra(\rho_0)$$
is well-defined and is an equivalence of groupoids. \label{prop:pe} \end{prop}

\proof To prove that this functor is well-defined as well as surjective on
isomorphism classes, we shall construct an inverse image $mon_b^{-1}(\rho)$ for
any $\rho\in\tra(\rho_0)$ and show that it does belong to $\tfa(\nabla_0)$.
(There does not exists a natural quasi-inverse functor to $(mon_b, \e_b)$,
though.) Our procedure is as follows: first we have the standard
construction of a flat bundle $W_A$ out of a representation $\pi_1(U, b)\to
G_A$. Then we consider the canonical extension $\bar{W}_A$ which is a
holomorphic bundle over $X$ with an $A$-action  such that all fibers are free
$A$-modules. Let us consider for any $l\in\pi_1(U, b)$ the following
diagram:  $$\begin{array}{ccc} A^N & \stackrel{\rho(l)}{\longrightarrow} & A^N
\\ \downarrow p & {} & \downarrow p \\ \C^N & \stackrel{\rho_0(l)}{\to} & \C^N
\end{array}.$$ This diagram is clearly commutative since $p(\rho)=\rho_0$. This
means that $\bar{W}_A\otimes_A \C$ is isomorphic to ${\bar W}_\C$.
Then Lemma \ref{lem:lNak} implies that $\bar{W}_A$ is (non-uniquely)
isomorphic to ${\bar W}_\C \otimes_\C A$. Besides, since $\rho(\g_i)$ is
conjugate in $G_A$ to $\rho_0(\g_i)$ for $1\le i\le r$ the eigenvalues (which a
priori belong to $A$) do not change and therefore the local system $\tag\otimes
A$ as well as $\tag$ splits locally into a direct sum of local system of rank
$1$ (over $A$).  Now one sees directly from the arguments used in \ref{lem:d'l}
that the residues of $mon_b^{-1}(\nabla)$ are the same as of $\nabla_0$. For the
rest of the proof one copies arguments from the proof of Proposition 6.3 in
\cite{GM2}. $\bigcirc$

\

Now we are ready to continue this section with a result similar in all ways to
Theorem 6.8 of \cite{GM2}. Let $L'$ denote the augmentation ideal
$$L':=Ker[B^{\fd}(\tag) \stackrel{k}{\to} \ag], $$ where $k$ is the
differential of the map $\e_b: \exp(B^0(\tag))\to G$ at the identity.
Analogously to the groupoid ${\tilde{\cal C}_A}$ one defines the groupoid
${\tilde{\cal C}_A(L')}$.

\begin{Th} The analytic germ of $\tr(\pi_1(U), G)_\N$ at $\rho_0$ pro-represents
the functor $$A\to Iso({\tilde{\cal C}_A(L')}).$$ \end{Th}

Now the proof of Theorem 1 in \cite{GM2} immediately implies that the functor
$$A\to Iso({\tilde{\cal C}_A(L')})$$ from the above Theorem is pro-represented
by the analytic germ of the quadratic cone consisting of all $u\in H^1(X,
j_*\tag )$ such that $[u,u]=0$. Therefore we conclude

\begin{Th} Let $\rho$ be a representation $\pi_1(U)\to G$ such that
$\rho(\g_i)\in\CC_i$. Then $\tr(\pi_1(U), G)_\N$ is quadratic at $\rho$.
\label{Th:mt} \end{Th}

\section{Moduli of representations of $\pi_1(U)$}
\setcounter{equation}{0}
In this section we deduce the consequences of results obtained earlier in this
paper to find out the geometric structure of the moduli space $\MN$. This space
was studied by J.-L. Brylinski and the author in \cite{BF}.  In that paper the
condition of vanishing of the second relative group cohomology group
$H^2(\pi_1(U), (\Gamma_i), \ag)$ was imposed to conclude that $\MN$ is smooth at
a point $[\rho_0]\in\MN$; here again $\ag$ is a $\pi_1(U)$-module via the
adjoint representation followed by $\rho_0$. Under the same assumption, the
tangent space $T_{[\rho]}\MN$ was identified as the first intersection
cohomology group (middle perversity) $IH^1(X, \tag)$. Other possible
descriptions of $T_{[\rho]}\MN$ include relative group cohomology group
$H^1(\pi_1(U), (\Gamma_i), \ag)$ the first cohomology group
$H^1(X, j_*\tag)$ of our controlling DGLA, and the first $L_2$ cohomology
group $H^1(U, \tag)_{(2)}$. In the same paper \cite{BF} it was
proven that if the smooth locus is reduced then it is symplectic and in the case
when $X$ is quasi-projective, then it actually is a \K manifold.

When we pass to the quotient $$\MN=\tr^{irr}(\pi_1(U), G)_\N/G,$$
due to the fact that we consider only irreducible representations and therefore
all $G$-orbits on $\tr^{irr}(\pi_1(U), G)_\N$ are compact and equidimensional,
we do not acquire additional singularities. Therefore Theorem \ref{Th:mt}
combined with Section 3 of \cite{BF} imply

\begin{Th} The singularities of the moduli space $\MN$ are quadratic.
\label {Th:modq} \end{Th}

Let us take a point $[\rho]\in\MN$; as was found in \cite{BF} for this point
to be smooth, vanishing of an infinite number of obstructions in $H^2(\pi_1(U),
(\Gamma_i), \ag)$ is sufficient. The first obstruction in this series is
given by the pairing
\begin{equation}
H^1(\pi_1(U), (\Gamma_i), \ag)\times H^1(\pi_1(U), (\Gamma_i),
 \ag) \to H^2(\pi_1(U), (\Gamma_i), \ag),
\label{eq:e37} \end{equation} defined as cup product in group
cohomology together with Lie bracket as coefficient pairing. Since we know now
that singularities of the space $\MN$ are quadratic, this obstruction is the
only one.

\begin{Th} A point $[\rho]\in\MN$ is a smooth point of $\MN$ if the pairing
\ref{eq:e37} is identically zero. \label{Th:pair} \end{Th}

The same conclusion about smoothness can be drawn if the pairing $$[,] :
H^1(X, j_*\tag )\times H^1(X, j_*\tag )\to H^2(X, j_*\tag )$$ is identically
zero.

\vskip 0.3in
Department of Mathematics \\ Penn State University \\ University Park, PA 16802
\\ foth@math.psu.edu

\

\noindent{\it AMS subj. class.}: \ \ primary 32G08;\ \ secondary 20L05, 17B70.

\end{document}